\documentclass[12pt,draft]{article}

\usepackage[latin1]{inputenc}
\usepackage[T1]{fontenc}
\usepackage{ae}
\usepackage{a4,amssymb,mathrsfs}
\def\Bbb{\mathbb}

\def\Zeta{{\rm Z}}
\def\zf{$\zeta$-function}
\def\fP{\mathfrak p}
\def\fA{\mathfrak a}

\def\cN{{\cal N}}
\def\vp{h}
\def\ph{\varphi}
\def\Tr{\mathop{\rm Tr}}
\def\Or{\mathop{\rm Or}\nolimits}
\def\OC{{\mathscr O}_C}
\def\cX{{\cal X}}
\def\cF{{\cal F}}
\def\ol{\overline}
\newtheorem{prop}{Proposition}
\newtheorem{lem}[prop]{Lemma}
\newtheorem{theorem}[prop]{Theorem}
\def\roep #1.{\medbreak\noindent{\it #1\/}.\enspace}
\def\endroep{\par\medbreak}
\def\qed{\hfill$\qquad\square$\par\medbreak}

\begin{document}
\title{\bf Poincaré series and zeta function for an irreducible
  plane curve singularity}
\author{Jan Stevens}
\date{}

\maketitle 
\begin{abstract}
The Poincaré series of an irreducible plane curve singularity
equals the \zf\ of its monodromy, by a result of Campillo,
Delgado and Gusein-Zade. We derive this fact from a formula of Ebeling
and Gusein-Zade relating the Poincaré
series of a quasi-homogeneous complete intersection 
singularity to  the Saito dual
of a product of \zf s.
\end{abstract}

Several cases are known where the $\zeta$-function of the monodromy of
an isolated hypersurface singularity is related to the
Poincaré series of its coordinate ring. 
The first instance of this phenomenon was observed by Campillo,
Delgado and Gusein-Zade \cite{CDG}: for an irreducible plane curve singularity
the $\zeta$-function of the monodromy equals its Poincaré series.
The proof is by comparing explicit formulas.
All attempts have failed at a more direct proof.
But maybe
this example is misleading. The results of Ebeling and
Gusein-Zade \cite{Eb, EGZ}  suggest a more indirect connection.
For quasi-homogeneous complete intersection 
singularities they show that the Poincaré
series, corrected by an orbit invariant, is related to the Saito dual
of a product of \zf s. The precise formulation will be given below.

The object of this note is to derive the original result of  \cite{CDG}
from the formula of \cite{EGZ}. 
The question
remains to give a direct proof of their formula and to explain the
meaning of the Saito dual.
The key observation, stressed by Teissier in \cite{Tei},  is that an
irreducible plane curve singularity is 
a (in some sense equisingular) deformation of a monomial curve with
the same semigroup, which is moreover a complete intersection. 
What is needed is a topological argument, which connects the \zf\ of
the plane curve singularity with the product of \zf s occurring in the
formula of \cite{EGZ}.

In the first section we give the definition of the monodromy \zf\
and motivate our conventions by recalling the analogy with number
theory.
Then we state the formula of Ebeling and Gusein-Zade. In the final
section we derive the main result.

\section{The zeta function of the monodromy}
There is some confusion in the literature about the definition of the
$\zeta$-function for singularities. 
We follow Milnor's book \cite{Mi2}.
The definition is based on an analogy with number theory.
For a number field $K$ the Dedekind $\zeta$-function is
$$
\zeta_K(s)= \prod_{\fP} \frac1{1-\frac1{\cN(\fP)^s}} =
\sum_{\fA} \frac1{\cN(\fA)^s}\;,
$$
where the product runs over all prime divisors $\fP$ of $K$
and the sum over all integral divisors, and $\cN(\fA)$ is the norm of
$\fA$.
For curves over finite fields a $\zeta$-function was defined by Emil
Artin in his thesis. In modified form (following F.K.~Schmidt) it is
given by exactly the same formula, where $K$ now stands for the function
field $K/{\Bbb F}_q$, and $\cN(\fP)$ is the number of elements in the
residue field, so equals $q^{\deg \fP}$.
A good reference is the second volume of Hasse's collected papers,
especially the paper in Italian \cite{Has}.

\roep Example. Let  $\zeta_0(s)$ be the $\zeta$-function of ${\Bbb
  P}^1$.
The number $N(n)$ of integral divisors of degree $n$ is the number of
homogeneous polynomials of degree $n$ with coefficients in ${\Bbb
  F}_q$ modulo scalars, so equals $\frac{q^{n+1}-1}{q-1}$. One 
finds
$$
\zeta_0(s)= \sum \frac{N(n)}{q^{ns}}
=\frac1{1-\frac q{q^s}} \frac1{1-\frac 1{q^s}} =\frac1{(1-qT)(1-T)}\;,
$$
where $T=q^{-s}$.
\endroep

André Weil  generalised the
function to schemes $X$ of finite type over ${\Bbb F}_q$, again by the
same product formula as a product over all closed points.
He conjectured that $\zeta_X(T)$ is a rational function of $T=q^{-s}$
(for more precise statements see \cite[Appendix C]{Har}).
This was known for curves.
If one sets $N_n=\#X({\Bbb F}_{q^n})$ then one
can compute that 
$$
\zeta_X(s) = \exp\;\sum_{n=1}^\infty \frac{N_n}n T^n\;. 
$$
The number $N_n$ is also the the number of fixed points of the $n$th
power of the Frobenius morphism $F$. The
Lefschetz fixed point formula in topology
computes the algebraic number  $\Lambda(\vp)$ of fixed points 
(and therefore $\#\mbox{Fix}(\vp)$ if all fixed points are
non-degenerate of index one) for a map $\vp\colon X\to X$ as
alternating sum of traces:
$$
\Lambda(\vp) = \sum (-1)^i \Tr \{\vp_*\mid H_i(X,{\Bbb
C})\}\;.
$$
An exercise in linear algebra 
gives
$$
\exp\;\sum_{n=1}^\infty \frac{\sum (-1)^i \Tr \{\vp_*^n\mid H_i\}}n T^n
=
\prod \det({\rm Id} - T\vp_*\mid H_i)^{(-1)^{i+1}} \;.
$$
Rationality of the \zf\ follows from  the existence of an appropriate
cohomology theory for varieties defined over fields of finite
characteristic.   

Milnor then used the same formula to define $\zeta$-functions in
topology
\cite{Mi1,Mi2}. 
Let $\vp\colon X\to X$  be a homeomorphism of a compact  
Euclidean Neighbourhood Retract in ${\Bbb R}^n$ (see
\cite[p.~81]{Dol}; for this type of spaces the Lefschetz-Hopf fixed
point theorem is proven in \cite[VII.6]{Dol}).
Sometimes the Lefschetz number $\Lambda(\vp)$ can be computed as the
Euler characteristic of the fixed point set  of $\vp$. This is the
case for any isometry of a compact Riemannian manifold
(cf.~\cite[Lemma 9.5]{Mi2}). It is true whenever the fixed point set
is a deformation retract of a $h$-invariant compact neighbourhood.
Let $\Lambda(\vp^k)$ be the Lefschetz number of the $k$th iterate of
$\vp$.
Define numbers $\chi_j$ by the recurrence $\Lambda(\vp^k)
= \sum_{j|k}\chi_j$. 
The number
$\chi_j$ is the Euler characteristic of the set of points of
primitive period $j$ (to have additivity of Euler characteristics also
for non closed sets we use Borel-Moore homology, see \cite[Sect.~91.1]{Ful}).
Then one can express the \zf\ of $\vp$ as 
$$
\zeta_\vp(T)=
\prod (1-T^j)^{-\chi_j/j}\;.
$$

We apply this to the monodromy of a function germ $f\colon X\to {\Bbb
C}$, defined on an equidimensional space. Consider a good
representative $f\colon \ol X\to D$ and let $F=f^{-1}(t)$ be the (possibly
singular)  Milnor fibre. In this situation one has a geometric
monodromy and a \zf\ $\zeta_f(T)$ of the monodromy.
To compute the \zf\ from the formula above one
needs a nice model for the monodromy. In the quasi-homogeneous case
it can be derived from the ${\Bbb C}^*$-action (cf.~\cite[Lemma 9.4]{Mi2}).
For a good embedded resolution one has A'Campo's
construction \cite{A'C}.

There is a fibre bundle of pairs $(\ol X,f^{-1}(t))$ over $D\setminus 0$ 
and a monodromy action
on the relative groups $H_*(\ol X,F)$. The corresponding \zf\
$\tilde\zeta_f(T)$ is related to the function $\zeta_f(T)$ of the monodromy
action on $H_*(F)$ by 
$$
\tilde\zeta_f(T)=(1-T)^{-1}(\zeta_f(T))^{-1}\;.
$$ 
It is also the inverse of the \zf\ on reduced homology.

\section{Zeta and Poincaré}
The relation between Poincaré series of quasi-homogeneous
singularities and $\zeta$-functions has been described by
Ebeling and Gusein-Zade \cite{EGZ}.

Let $A=\oplus A_k$ be a graded ring with each $A_k$ finite
dimensional. Its Poincaré series is 
$P_A(t)=\sum (\dim A_k) T^k$. For a general ring with a filtration
the Poincaré series is by definition the  Poincaré series
of its associated graded ring. The Poincaré series of a singularity
is the   Poincaré series of its local ring.

Let $X$ be a weighted homogeneous
complete intersection in ${\Bbb C}^n$ given by the ideal
$(f_1,\dots,f_k)$ with $\deg f_i=d_i$ and with wt $z_i=q_i$.
Set $X^{(j)}=V(f_1,\dots,f_j)$ for $j=0,\dots,k$ (so $X^{(0)}={\Bbb
  C}^n$ and
$X^{(k)}=X$).
Consider the function $f_j$ on the space  $X^{(j-1)}$ and let
$F^{(j)}$ be its Milnor fibre:
 $F^{(j)}= f_j^{-1}(1) \cap X^{(j-1)}$.
As  monodromy transformation one can take 
$(z_1,\dots,z_k)\mapsto (e^{2q_1\pi/d_j}z_1,\dots,e^{2q_n\pi/d_j}z_n )$.
The function $\tilde\zeta_j(T)$ is the $\zeta$-function
of the monodromy action on $H^*(X^{(j-1)},F^{(j)})$.
Let  $Y = (X\setminus 0)/{\Bbb C}^*$ be the orbit space
and let $Y_m$ be the set of
orbits with isotropy group ${\Bbb Z}/ m {\Bbb Z}$.
One defines the orbit invariant 
$$
\Or_X(T) = \prod_{m\geq1}(1-T^m)^{\chi(Y_m)}\;.
$$

Following Kyoji Saito \cite{S1,S2} one defines the Saito dual (of
level $d$) of a 
rational function 
$\vp(T)=\prod_{l|d}(1-T^l)^{\alpha_l}$ to be 
$$
\vp^{*_d}(T)=\prod\nolimits_{m|d}(1-T^{d/m})^{-\alpha_{m}}\;, 
$$ 
where the notation emphasises the dependence on the level $d$.

We can now state the result of \cite{EGZ}:
\begin{theorem} \label{propegz}
$$
P_X(T) \Or_X(T)=\prod_{j=1}^k (\tilde\zeta_j)^{*_{d_j}}(T) \;.
$$
\end{theorem}

\section{Irreducible plane curve singularities}
Let $(C,0)$ be an irreducible plane curve singularity. The $t$-adic
valuation $\nu$ on ${\Bbb C}\{t\}$, the normalisation of $\OC$, 
induces a filtration on $\OC$. Let gr$\;\OC$ be the associated graded
algebra. Monique Lejeune-Jalabert (see Teissier's Appendix in
\cite{Tei})  showed
that it is a graded algebra isomorphic to ${\Bbb C}[C^\Gamma] = {\Bbb
  C}[t^h\colon h\in\Gamma]$, the ring of the monomial curve $C^\Gamma$
with the same semigroup $\Gamma=\langle \beta_0,\dots,\beta_g\rangle$
as $C$. 

\begin{prop}
The curve $C^\Gamma\subset {\Bbb C}^{g+1}$ is a complete
intersection. 
\end{prop}
 
A proof can be found in Teissier's paper. We describe the equations,
using the structure of the semigroup $\Gamma=\langle
\beta_0,\dots,\beta_g\rangle$. One sets $e_0=\beta_0=n$,
$e_i=\gcd(e_{i-1},\beta_i)$, $e_{i-1}=n_ie_i$.
Then $\beta_0=n_1\cdots
n_g$, $n_i\beta_i<\beta_{i+1}$ and $n_i\beta_i \in \langle
\beta_0,\dots,\beta_{i_1}\rangle$, so $n_i\beta_i=\sum
l_{ij}\beta_j$.  Write $z_i=t^{\beta_i}$. The curve $C^\Gamma$ has equations
$$
f_i=z_i^{n_i}-\prod z_j^{l_{ij}}\;.
$$
According to
Teissier [loc.cit.], the curve $C$ occurs in the versal deformation of
$C^\Gamma$ (in the part with positive weights.) 
In fact, each ring containing a field is a deformation of its
associated graded w.r.t.~any filtration. One gets also an easy way to
write down equations of a plane curve with the same semigroup by
taking the special 
deformation $\ph_i=f_i +\lambda_i z_{i+1}$ with $\lambda_i$ a suitable
power of a deformation variable $s$. The equations $\ph_1,\dots,\ph_{g-1}$
define a smooth space on which the function $f_g$ gives a plane curve
equisingular with the original $C$.

\roep Example.
Let $\Gamma=\langle 4,6,13 \rangle$ and set
$(x,y,z)=(t^4,t^6,t^{13})$.
Equations are $y^2-x^3=0$ and $z^2-x^5y=0$.
We deform the first equation into $y^2-x^3=sz$. For $s\neq0$ we can
eliminate $z$ and find the plane curve $(y^2-x^3)^2-s^2x^5y=0$.
\endroep

The result of Campillo, Delgado and Gusein-Zade \cite{CDG} is:
\begin{theorem}
\label{propcdg}
For an irreducible plane curve singularity the $\zeta$-function of its
monodromy equals its Poincaré series. 
\end{theorem}

We shall derive it from the formula of Ebeling and Gusein-Zade
\cite{EGZ} (Theorem \ref{propegz}).
The Poincaré series in question is the same as that of $C^\Gamma$.
The orbit invariant reduces to $\Or_{C^\Gamma}(T)=(1-T)$.
To distinguish between the \zf s occurring in the formula and
monodromy \zf s of plane curve singularities we will denote the latter
by capital $\Zeta$.
Theorem \ref{propcdg} is equivalent to the following statement:
\begin{prop}
\label{propjan}
With the notations introduced above
$$
(\widetilde\Zeta_C(T))^{-1} =\prod_{j=1}^g (\tilde\zeta_j)^{*_{d_j}}(T) \;.
$$
\end{prop}

The proof will be by induction on $g$. Therefore we denote $C$ by
$C_g$, its Milnor fibre by $F_g$ and \zf\ by $Z_g(T)$.
Let $C_j$ be the $j$th approximate curve to $C_g$ (with Milnor fibre
$F_{j}$ and \zf\ $Z_{j}(T)$). It has $j$
Puiseux pairs, which equal the first $j$ Puiseux pairs of $C_g$.
The semigroup $\Gamma_j$ of $C_j$ has $j+1$ generators. In particular,
$\Gamma_{g-1}= \langle
\beta_0/n_g,\dots,\beta_{g-1}/n_g\rangle$. 
We denote the corresponding monomial curve of
embedding dimension $j+1$ by $C^{(j)}$. 
It has the same
equations as the space $X^{(j)}$, since the variables $z_{j+1}$, \dots,
$z_g$ do not occur in the first $j$ equations. Therefore $X^{(j)}=
C^{(j)} \times {\Bbb C}^{g-j}$. 
The monodromy action on the second factor is
$(z_{j+1},\dots,z_g)\mapsto (e^{2q_{j+1}\pi/d_j}z_{j+1},
     \dots,e^{2q_n\pi/d_j}z_g )$.
The contribution to the Euler characteristic of the fixed point set,
which comes from  $C^{(j)} \times ({\Bbb C}^{g-j}\setminus 0)$, is
zero.
For the \zf\ we can therefore forget about the factor.
In particular, the curves $C^{(g)}$ and
$C^{(g-1)}$ give rise to the same \zf s $\tilde\zeta_1$, \dots,
$\tilde\zeta_{g-1}$. But the degrees of the equations are different:
if the equation $f_j$ of $C^{(g)}$ has degree $d_j$, then it
has degree $d_j/n_g$ considered as equation of $C^{(g-1)}$.
In the Ebeling--Gusein-Zade formula (Theorem \ref{propegz}) the
Saito dual is therefore taken at a different level. One has
$$ 
(\tilde\zeta_j)^{*_{d_j}}(T)=(\tilde\zeta_j)^{*_{d_j/n_g}}(T^{n_g})\;.
$$
This said, Proposition \ref{propjan} follows by induction from the
following two lemmas.

\begin{lem}
\label{lem1}
$\quad(\tilde\zeta_g)^{*_{d_g}}(T)=(\tilde\zeta_g(T))^{-1}\;.$
\end{lem}
\begin{lem}
\label{lem2}
$\quad\widetilde \Zeta_g(T) 
   = \tilde\zeta_g(T)\cdot \widetilde \Zeta_{g-1}(T^{n_g})\;.$
\end{lem}

\roep Proof of Lemma {\rm\ref{lem1}}.
The space $X^{(g-1)}$ is homeomorphic to its normalisation, which is
smooth. In coordinates $(\tau,z_g)$ on the normalisation the function
$f_g$ is given by $z_g^{n_g}-\tau^{\beta_g}$; it has degree
$d_g=n_g\beta_g$. One easily computes that
$$\tilde
\zeta_g(T)=\frac{(1-T^{\beta_g})(1-T^{n_g})}{(1-T^{d_g})(1-T)}
\;,
$$
so $\tilde\zeta_g$ and its inverse are dual of level $d_g=n_g\beta_g$.
\qed

\roep Proof of Lemma {\rm\ref{lem2}}.
We consider the degeneration of $C_g$ to the monomial curve
with deformation parameter $s$. We have a complete intersection in
$({\Bbb C}^{g+1}\times{\Bbb C},0)$ defined by functions
$\ph_1(z,s)$, \dots, $\ph_g(z,s)$. Consider the function $\Phi_g=(\ph_g,s)$ on
${\cX}^{(g-1)}=V(\ph_1,\dots,\ph_{g-1})$. For fixed $s\neq0$ it
defines the curve $C_g$ and for $s=0$ the monomial curve $C^{(g)}$.
The Milnor fibre $F_g$ is $\Phi^{-1}(t,s)$ (for $t,s\neq0$),
while $\Phi^{-1}(t,0)$ (for $t\neq0$) gives the singular Milnor fibre
$F^{(g)}$ with $n_g$ singular points of type $C^{(g-1)}$ (the
transversal singularity of $X^{(g-1)}$). The intersection of a small
ball around such a singular point with the Milnor fibre $F_g$ (for
small $s$) is a Milnor fibre $F_{g-1}$ of the plane curve singularity
$C_{g-1}$. 

We are interested in the $t$-monodromy. For $\tilde\zeta_g(T)$
we have to consider the pair $(X^{(g-1)},F^{(g)})$.
The fibre $F^{(g)}$ is a retract of $\cF^{(g)}=\Phi^{-1}(t\times
D_\delta)$, where $D_\delta$ is a small disc around the origin. In the
same way we construct a thickening $\cX^{(g-1)}$ of $X^{(g-1)}$.
Then $H_*(\cX^{(g-1)},\cF^{(g)})=H_*(X^{(g-1)},F^{(g)})$.
The Milnor fibre $F_g$ is now a subspace of $\cF^{(g)}$, so it makes
sense to consider the long exact sequence of the triple
$(\cX^{(g-1)},\cF^{(g)},F_g)$:
$$
\cdots \to
H_q(\cF^{(g)},F_g)\longrightarrow
H_q(\cX^{(g-1)},F_g)\longrightarrow
H_q(\cX^{(g-1)},\cF^{(g)})\to
\cdots
$$
As $\cX^{(g-1)}$ is also contractible one has 
$H_*(\cX^{(g-1)},F_g)\cong H_*(X_g,F_g)$, where $X_g$ is the the
intersection of $\cX^{(g-1)}$ with the hyperplane $\{s=s_0\}$ for
$s_0$ the same value used to define the Milnor fibre $F_g$, and the
isomorphism is compatible with the monodromy action.
Outside small neighbourhoods of the singular points of $F^{(g)}$ the
thickened fibre $\cF^{(g)}$ is a trivial fibre bundle over $D_\delta$
so by deformation retraction and excision
$H_*(\cF^{(g)},F_g)\cong \oplus_{i=1}^{n_g} H^*(X_{g-1,i}, F_{g-1,i})$
where $X_{g-1,i}$ is a small Milnor ball around the $i$th singularity. 
The monodromy permutes the Milnor balls, and the action after $n_g$ steps 
on $H^*(X_{g-1,i}, F_{g-1,i})$ is the monodromy of the curve
singularity $C_{g-1}$. Therefore the \zf\ of the monodromy on
$\oplus_{i=1}^{n_g} H_*(X_{g-1,i}, F_{g-1,i})$ is 
$\widetilde\Zeta_{g-1}(T^{n_g})$.
The lemma now follows from the exact sequence
$$
\cdots \to
\oplus_{i=1}^{n_g} H_q(X_{g-1,i}, F_{g-1,i})\longrightarrow
H_q(X_g,F_g)\longrightarrow
H_q(X^{(g-1)},F^{(g)})\to
\cdots
$$
\qed

\vfill

\noindent Matematik\\
  Chalmers tekniska h\"ogskola och G\"oteborgs universitet, \\
  SE 412 96 G\"oteborg, Sweden\\
  e-mail: \texttt{stevens@math.chalmers.se}

\end{document}